\documentstyle[12pt]{amsart}
\newtheorem{th}{Theorem}[section]
\newtheorem{crl}[th]{Corollary}
\newtheorem{prp}[th]{Proposition}
\newtheorem{lm}[th]{Lemma}

\newcommand{\der }{\partial }
\newcommand{\eps}{\epsilon}
\begin{document}
\title{Minimal polynomial identities for right-symmetric
algebras}
\author {Askar Dzhumadil'daev}
\address
{Institute of Mathematics, Pushkin str.125, Almaty, 480021,
Kazakhstan}
\email{askar@@itpm.sci.kz}
\maketitle
\begin{abstract}
An algebra $A$ with multiplication 
$A\times A\rightarrow A, (a,b)\mapsto a\circ b,$ 
is called right-symmetric, if $a\circ(b\circ c)-(a\circ b)\circ c=
a\circ (c\circ b)-(a\circ c)\circ b,$ for any $a,b,c\in A.$ 
The multiplication of right-symmetric Witt algebras 
$W_n=\{u\der_i: u\in U, U={\cal K}[x_1^{\pm 1},\ldots,x_n^{\pm 1}],$ 
or $={\cal K}[x_1,\ldots,x_n], i=1,\ldots,n\}, \;
p=0,$ or $W_n({\bf m)}=\{u\der_i: u\in U, U=O_n({\bf m})\}, \; p>0,$ 
are given by $u\der_i\circ v\der_j=v\der_j(u)\der_i.$
An analogue of the Amitsur-Levitzki theorem for right-symmetric 
Witt algebras is established. Right-symmetric Witt algebras of rank $n$ 
satisfy the standard right-symmetric identity of degree $2n+1:$ 
 $\sum_{\sigma\in 
Sym_{2n}}sign(\sigma)a_{\sigma(1)}\circ(a_{\sigma(2)}\circ
\cdots(a_{\sigma(2n)}\circ a_{2n+1})\cdots)=0.$ The minimal degree for
left polynomial identities of $W_n^{rsym}, W_n^{+rsym}, p=0,$ is $2n+1.$ 
The minimal degree of multilinear left polynomial identity of $W_n({\bf m}), p>0,$ 
is also $2n+1.$ All left polynomial (also multilinear, if 
$p>0$) identities of right-symmetric Witt algebras of  minimal degree are linear
combinations of left polynomials obtained from standard ones by 
permutations of arguments.
\end{abstract}

\section{Introduction}
According to the Amitsur-Levitzki theorem~\cite{Amitsur} the matrix algebra $Mat_n$ satisfies a standard polynomial identity 
of degree $2n:$ 
$$\sum_{\sigma\in Sym_{2n}}sign\,\sigma\,a_{\sigma(1)}\circ\cdots
a_{\sigma(2n)}=0,$$
where $a\circ b$ is a usual matrix multiplication. Moreover, $Mat_n$ has no 
polynomial identity of degree less that $2n.$ For details on polynomial 
identities of associative algebras see for example,~\cite{Procesi}.

An algebra $W_1=\{e_i: e_i\circ e_j=(i+1)e_{i+j}, i,j\in {\bf Z}\}$ is
right-symmetric. Since 
its Lie algebra is isomorphic to a Witt algebra $W_1=\{e_i:
[e_i,e_j]=(j-i)e_{i+j}\}$ 
 we call it as right-symmetric Witt algebra of rank 1 and denote by
$W_1^{rsym}.$ This algebra satisfies a right-symmetric identity
$$a\circ(b\circ c)-(a\circ b)\circ c=a\circ (c\circ b)-(a\circ c)\circ b$$
and a left-commutativity identity
\begin{equation}
a\circ(b\circ c)=b\circ(a\circ c). \label{Novikov}
\end{equation}
Such algebras are called {\it Novikov} \cite{Osborn1}, \cite{Osborn2},
\cite{Osborn3}.

There is a generalisation of the Witt algebra to the many valuables case.
Let $U$ be an associative commutative algebra with a set of commuting derivations 
${\cal D}=\{\der_i: i=1,\ldots,n\}.$ For any $u\in U,$ an endomorphism $u\der_i: U\rightarrow U,$ 
such that $(u\der_i)(v)=u\der_i(v),$ is a derivation of $U.$ 
Denote by $U{\cal D}$ a space of derivations $\sum_{i=1}^nu_i\der_i.$
Endow this space by multiplication
$$u\der_i\circ v\der_j=v\der_j(u)\der_i.$$
We obtain a right-symmetric algebra $U{\cal D}.$ This algebra is called a 
{\it a right-symmetric Witt algebra generated on $U$ and ${\cal D}$.} 

In our paper, $U$ is ${\cal K}[x_1,\ldots,x_n],$ or Laurent polynomial algebra 
${\cal K}[x_1,^{\pm 1},\ldots,x_n^{\pm 1}],$ or a divided power algebra $O_n({\bf m})=\{x^{\alpha}:
x^{\alpha}x^{\beta}={\alpha+\beta\choose\alpha}x^{\alpha+\beta}\},$ if
the characteristic of ${\cal K}$ is $p>0.$ As a Lie algebra the Witt algebra 
of the rank $n$ is defined as a Lie algebra of derivations of $U.$ 
The multiplication $u\der_i\circ v\der_j=v\der_j(u)\der_i$ 
satisfies the right-symmetry identity. Obtained 
right-symmetric Witt algebras of rank $n,$ are denoted by $W_n^{rsym}$ or 
$W_n^{+rsym}$ or $W_n{(\bf m})^{rsym},$ depending on $U={\cal K}[x_1^{\pm 1},\ldots,
x_n^{\pm 1}],$ or ${\cal K}[x_1,\ldots,x_n],$ or $O_n({\bf m}).$
It is easy to notice that the right-symmetric Witt algebras of rank $n$ do
{\it not} satisfy the left-commutativity identity if $n>1.$ 

We are interested in the analogues of the left-symmetric identities for 
the case of many valuables. We suggest two ways to solve this problem. 

In the first way we endow a vector space of the Witt algebra with two 
multiplications: the multiplication $(a,b)\mapsto a\circ b,$ 
mentioned above and the second multiplication defined by
$u\der_i\ast v\der_j=\der_i(u)v\der_j.$ We obtain an algebra with the following 
identities
$$a\circ (b\circ c)-(a\circ b)\circ c-a\circ(c\circ b)+(a\circ c)\circ b=0,$$
\begin{equation}
a\ast (b\ast c)-b\ast(a\ast c)=0,\label{Novikov1}
\end{equation}
$$a\circ (b\ast c)-b\ast(a\circ c)=0,$$
$$(a\ast b-b\ast a- a\circ b+b\circ a)\ast c=0,$$
$$(a\circ b- b\circ a)\ast c+a\ast(c\circ b)-(a\ast c)\circ b
-b\ast(c\circ a)+(b\ast c)\circ a=0.$$
Notice that (\ref{Novikov1}), for the multiplication $\ast,$ is similar 
to the identity (\ref{Novikov}) for the multiplication $\circ.$

In the second way we save right-symmetric multiplication $(a,b)\mapsto a\circ b,$ 
and try to construct an identity for its left multiplication operators. 
For $a\in A,$ denote by $r_a$ and $l_a$ operators of right and left multiplications 
on $A:$ $br_a=b\circ a,\; bl_a=a\circ b.$ In terms of right and left multiplication 
operators the right-symmetry identity is equivalent to the following conditions
$$[r_a,r_b]=r_{[a,b]},\;\; [r_a,l_b]=l_al_b-l_{b\circ a}, \; \forall a,b\in A. $$

Our main result is the following. For right-symmetric Witt algebras of any 
rank $n$ for left multiplication operators the following {\it standard polynomial 
identity of degree $2n$} holds:
$$\sum_{\sigma\in Sym_{2n}}sign\,\sigma\,l_{a_{\sigma(1)}}\cdots l_{a_
{\sigma(2n)}}=0,$$
or in terms of multiplication $\circ,$ the following {\it left polynomial identity 
of degree $2n+1$} is valid:
$$\sum_{\sigma\in Sym_{2n}}sign\,\sigma\,a_{\sigma(2n)}\circ (a_{\sigma(2n-1)}\circ\cdots 
\circ(a_{\sigma(2)}\circ(a_{\sigma(1)}\circ a_0))\cdots)=0.$$
We prove that in the space of left polynomial identities the degree $2n+1$ 
for right-symmetric Witt algebras of rank $n$ is {\it minimal} (in the case of $p>0$ 
we suppose also that polynomials are multilinear). We also prove that the  
left polynomial identities of right-symmetric Witt algebras of rank $n$ of 
minimal degree {\it can be obtained from standard polynomials by permutation 
and linear combination operations.} 

Notice that, for $n=1$ the identity $s_{2}^{rsym}(a_0,a_1,a_2)=0,$ 
coincides with the Novikov identity 
$$a_0\circ(a_1\circ a_2)-a_1\circ(a_0\circ a_2)=0.$$ 
So, we can consider right-symmetric algebras that satisfy 
standard right-symmetric identity
\begin{equation}\sum_{\sigma\in Sym{2n}}sign\,{\sigma}\,a_{\sigma(1)}\circ
a_{\sigma(2)}\circ
\cdots a_{\sigma(2n)}\circ a_{0}=0,\label{ask}
\end{equation}
as a generalisation of Novikov algebras. This class of algebras includes
Witt algebras in the many valuables case. 

In our proof we use some properties of Laurent or divided
power polynomials. The identity~(\ref{Novikov}) is true for any associative commutative 
algebra $U$ and for a set ${\cal D}$ with one derivation $\der_1.$  We believe that the 
identity~(\ref{ask}) holds for any associative commutative algebra $U$ with 
a set of commuting derivations ${\cal D}$ with $n$ derivations $\der_i, i=1,\ldots,n.$ 

{\bf Conjecture.} {\it Any simple right-symmetric algebra over an algebraically closed field 
of characteristic $p=0$ or $p>2n+1$ with minimal left 
polynomial identity $s_{2n}^{rsym}$ is isomorphic to one of the following algebras:
\begin{itemize}
\item a Witt algebra generated by some associative commutative algebra $U$ with 
a set of linear independent commuting derivations ${\cal D}$ with $n$ derivations, 
\item or their deformations. 
\end{itemize}}
About deformations of right-symmetric algebras and the description of local 
deformations of $A=W_n^{rsym}, W_n^{+rsym},$ if $p=0,$ or $W_n({\bf m}),$ if $p>0,$ 
see~\cite{DzhumaKos}, \cite{Gerst}.  As an example let us give some right-symmetric deformations 
of $A=W_1=\{a=u\der: u\in {\cal K}[x^{\pm 1}]\}.$ 

The space of local deformations of $A$ is 4-dimensional and generated by classes 
of the following right-symmetric 2-cocycles
$$\psi^1(u\der,v\der)= x^{-1}uv\der,$$
$$\psi^2(u\der,v\der)= x^{-1}\der(u)v\der,$$
$$\psi^3(u\der,v\der)= (u-x\der(u))\der(v)\der,$$
$$\psi^4(u\der,v\der)= \der(u)\der(v)\der.$$
If $\Psi_1=\sum_{k=1}^4\varepsilon_k\psi^k$ is a 4-parametrical 
local deformation of $A,$ is it possible to construct prolongations $\Psi_l=
\sum_{|{\bf i}|=l}\eps^{\bf i}\psi^{\bf i},$ where ${\bf i}=(i_1,i_2,i_3,i_4), 
|{\bf i}|=i_1+i_2+i_3+i_4?$ In other words, is it possible to find 
$\psi^{\bf i}\in C^2_{rsym}(A,A),$ such that a new multiplication
$$a\circ_{\bf \varepsilon}b=a\circ b+\sum_l\Psi_l,$$
will be a right-symmetric multiplication over a field 
${\cal K}(({\bf \varepsilon}))?$ The answer is: $\Psi_1$ can be prolongated to a 
global deformation if and only if $\varepsilon_1\varepsilon_4+
\varepsilon_2\varepsilon_3=0.$ We give prolongation formulas for some 
special cases.

The local deformation $\varepsilon_1\psi^1+\varepsilon_2\psi^2$ 
of $W_1$ has a trivial prolongation:
$$(u\der,v\der)\mapsto
\der(u)v+\varepsilon_1x^{-1}uv+\varepsilon_2x^{-1}\der(u)v,$$ 
is a right-symmetric multiplication. This algebra was obtained by Osborn 
\cite{Osborn3}. His results confirm our conjecture for the case of $n=1.$ 
Notice that, cocycles $\psi^3,$ $\psi^4$  do not 
satisfy left-commutativity identity. They are not Novikov cocycles \cite{Nov}.
Each of these cocycles have the following prolongations:
\begin{equation}
(a,b)\mapsto \der(a)b+[x\der,a](\sum_{i}\varepsilon_3^i(-1)^ix^i\der^i/
\{(\varepsilon_3+1)\cdots(i\varepsilon_3+1)\})(b),\label{eps3}
\end{equation}
$$
(a,b)\mapsto \der(a)\sum_{i}\varepsilon_4^i\der^i(b).
$$
In the case of (\ref{eps3}) we should change expressions like 
$(i\varepsilon_4+1)^{-1}$ to formal series
$1-i\varepsilon_3+i^2\varepsilon_3^2-i^3\varepsilon_3^3+\cdots .$ 
Then we obtain a formal power serie
$$(a,b)\mapsto \der(a)b$$
$$-\varepsilon_3[x\der,a]\der(b)$$
$$+\varepsilon_3^2[x\der,a](x\der(b)+x^2\der^2)(b)$$
$$-\varepsilon_3^3[x\der,a](x\der+3x^2\der^2+x^3\der^3)(b)$$
$$+\varepsilon_3^4[x\der,a](x\der+7x^2\der^2+6x^3\der^3+x^4\der^4)(b)$$
$$-\varepsilon_3^5[x\der,a](x\der+63x^2\der+25x^3\der^3+
10x^4\der^4+x^5\der^5)(b)+\cdots.$$
This is one of the prolongations of the local deformation $-\varepsilon_3[x\der,a]
\der(b).$ It will be interesting to construct prolongation formulas for a 
linear combination of cocycles and find polynomial identities of obtained algebras. 
It is also interesting to find right polynomial identities of right-symmetric algebras.
Right multiplication operators satisfy Lie algebraic conditions and in 
this case one can expect Lie algebraic difficulties (see~\cite{Razmyslov}, \cite{Razmyslov1}). 
Let us mention that an identity of degree 5 for Lie algebra $W_1^+$ 
is true also for right multiplication operators:
$$\sum_{\sigma\in Sym_4}sign\,\sigma\,r_{a_{\sigma(1)}}r_{a_{\sigma(2)}}
r_{a_{\sigma(3)}}r_{a_{\sigma(4)}}=0.$$
Moreover, for right-symmetric Witt algebra $U{\cal D},$ ${\cal D}=\{\der\},$ the 
following right polynomial identity of degree 3 takes place
$$\sum_{\sigma\in Sym_3}sign\,\sigma (a_{\sigma(1)}\circ a_{\sigma(2)})\circ 
a_{\sigma(3)}=0.$$
For a right-symmetric algebra $U{\cal D},$ where ${\cal D}=\{\der_1,\der_2\},$
the following right polynomial identity of degree 7 is true
$$\sum_{\sigma\in Sym_7}sign\,\sigma (((((a_{\sigma(1)}\circ a_{\sigma(2)})\circ 
a_{\sigma(3)})\circ a_{\sigma(4)})\circ a_{\sigma(5)})\circ a_{\sigma(6)})\circ 
a_{\sigma(7)} =0.$$

\section{Right-symmetric algebras}
Let $A$ be an algebra with multiplication 
$A\times A\rightarrow A, (a,b)\mapsto a\circ b.$ Let $(a,b,c)=a\circ(b\circ c)-
(a\circ b)\circ c$ be an associator of elements $a,b,c\in A.$ Associative algebras
are defined by condition $(a,b,c)=0,$ for any $a,b,c\in A.$ Right-symmetric 
algebras are defined by identity
$$(a,b,c)=(a,c,b),$$
i.e., by identity
$$a\circ(b\circ c)-(a\circ b)\circ c=a\circ (c\circ b)-(a\circ c)\circ b.$$
The left-symmetric identity is
$$(a,b,c)=(b,a,c).$$
There is one-to-one correspondence between right-symmetric and left-symmetric algebras.
Namely, if $(a,b)\rightarrow a\circ b,$ is right(left)-symmetric, then a new 
multiplication $(a,b)\mapsto b\circ a,$ is left(right)-symmetric. In our paper 
left-symmetric algebras are not considered. Right-symmetric algebras are
called sometimes as Vinberg-Kozsul algebras \cite{Vinberg}, \cite{Kozshul}.
Right-symmetric algebra $A$ is Lie-admissible, i.e., under commutator 
$[a,b]=a\circ b-b\circ a,$ we obtain a Lie algebra.

Any associative algebra is right-symmetric.
In a such cases, we will use notations like $A^{ass},$ if we consider
$A$ as an associative algebra and $A^{rsym}$ if we consider $A$ as 
right-symmetric algebra. Similarly, for a right-symmetric algebra
$A$ the notation $A^{rsym}$ will mean that we use only
right-symmetric structure on $A$ and $A^{lie}$ stands for
a Lie algebra structure under commutator $(a,b)\mapsto [a,b].$

In terms of operators of right multiplication $r_a$ and 
left multiplication $l_a$ of the algebra $A,$ 
$$ar_b=a\circ b, \; al_b=b\circ a,$$
right-symmetry identities are equivalent to the following conditions:
$$[r_a,r_b]-r_{[a,b]}=0,$$
$$[r_a,l_b]-l_al_b+l_{b\circ a}=0.$$
Let $A_{r}=\{a_r: a\in A\}$ and $A_l=\{a_l: a\in A\}$ be two copies of $A,$
and ${\cal A}=A_r\oplus A_{l}$ be their direct sum.
Let $T({\cal A})={\cal K}\oplus {\cal A}\oplus {\cal A}\otimes {\cal A}\oplus 
\cdots$ be the tensor algebra of ${\cal A}.$ 
A universal enveloping algebra of $A,$ denoted by $U(A),$ is defined 
as a factor-algebra of $T(A_r\oplus A_l)$ over an ideal generated 
by $[a_r,b_r]-{[a,b]}_r, [a_r,b_l]-a_lb_l+{(b\circ a)}_l.$ Denote elements of 
$U(A)$ corresponding to $a_r,a_l$ by $r_a,l_a.$ 

Let $B$ be a subalgebra of $A.$ Let 
$$Z^{l.ass}_A(B)=\{a\in A: (a,b_1,b_2)=0,\;\forall b_1,b_2\in B\},$$
$$Z^{r.ass}_A(B)=\{a\in A: (b_1,b_2,a)=0,\;\forall b_1,b_2\in B\},$$
be {\it left and right associative centralisers of $B$ in $A.$} Let
$$N^{l.ass}_A(B)=\{a\in A: (a,b_1,b_2)\in B, \forall b_1,b_2\in B\},$$
$$N^{r.ass}_A(B)=\{a\in A: (b_1,b_2,a)\in B, \forall b_1,b_2\in B\},$$
be {\it left and right normalisers of $B$ in $A.$} It is clear that
$$Z^{l.ass}_A(B)\subseteq N^{l.ass}_A(B),$$
$$Z^{r.ass}_A(B)\subseteq N^{r.ass}_A(B).$$
Let
$$Z^{left}_A(B)=\{a\in A: a\circ b=0,\;\forall b\in B\},$$
$$Z^{right}_A(B)=\{a\in A: b\circ a=0, \; \forall b\in B\},$$
be {\it left and right centralisers of $B$ in $A$} and 
$$N^{left}_A(B)=\{a\in A: a\circ b\in B,\;\forall b\in B\},$$
$$N^{right}_A(B)=\{a\in A: b\circ a\in B, \; \forall b\in B\},$$
be {\it left and right normalisers of $B$ in $A.$} We have
$$Z^{left}_A(B)\subseteq Z^{l.ass}_A(B),$$
$$Z^{right}_A(B)\subseteq Z^{r.ass}_A(B),$$
$$N^{left}_A(B)\subseteq N^{l.ass}_A(B),$$
$$N^{right}_A(B)\subseteq N^{r.ass}_A(B).$$

For the left cases and if $A=B$ we reduce these denotions: 
$Z(A)=Z^{left}_A(A),$ $ N(B)=N^{left}_A(B),$
$Z^{right}(A)=Z^{right}_A(A),$ 
$N^{right}(A)=N^{right}_A(A),$
$Z^{l.ass}(A)=Z^{l.ass}_A(A), N^{l.ass}(A)=N^{l.ass}_A(A),$
$Z^{r.ass}(A)=Z^{r.ass}_A(A), N^{r.ass}(A)=N^{r.ass}_A(A).$

We call $Z(A)$ and $Z^{right}(A)$ as left and right centers of $A.$
We call also $Z^{l.ass}(A)$ and $Z^{r.sym}(A)$ as left and right 
associative centers of $A.$ 

Left (right) associative centers are close under multiplication $\circ.$
To see this let us consider for simplicity the case of left associative centers.
Suppose that $X,Y\in Z^{l.ass}(A).$ Then according to the right-symmetric identity
$$(X\circ Y)\circ (a\circ b)=$$
$$X\circ (Y\circ (a\circ b))=
X\circ ((Y\circ a)\circ b)=(X\circ (Y\circ a))\circ b=$$
$$((X\circ Y)\circ a)\circ b.$$
So, $X\circ Y\in Z^{l.ass}(A),$ and $Z^{left}(A)$ is a subalgebra of $A.$

Notice that $Z(A)$ and $N(Z(A))$ are also subalgebras of $A:$
$$(z_1\circ z_2)\circ a=z_1\circ(z_2\circ a)-(z_1,a,z_2)=0,$$
$$((n_1\circ n_2)\circ z)\circ a=(n_1\circ (n_2\circ z))\circ a
-(n_1,z,n_2)\circ a=0,$$
for any $a\in A, z_1,z_2\in Z(A), n_1,n_2\in N(Z(A)).$ The same is true 
for $Z^{right}(A)$ and $N^{right}(Z^{right}(A)).$ 

\begin{prp} 
If $z\in Z(A),$ then $r_z$ is a derivation of $A.$
\end{prp}
{\bf Proof.} Since $z\circ b=0,$ we have $a\circ (z\circ b)=0.$
According right-symmetric identity
$$(a\circ b)\circ z= a\circ (b\circ z)+(a\circ z)\circ b,$$
for any $a,b\in A. \bullet$

\begin{prp}
For any $z\in Z(A), a\in N(Z(A)),$ and 
for any $b\in A,$ 
$$a\circ (b\circ z)=(a\circ b)\circ z.$$
\end{prp}
{\bf Proof.} Let $z\in Z(A).$ Then $z\circ b=0,$ and 
$a\circ(z\circ b)=0.$ Let $a\in N(Z(A)).$ 
Then $(a\circ z)\circ b=0.$ So, 
$$a\circ (b\circ z)-(a\circ b)\circ z= 
a\circ (z\circ b)-(a\circ z)\circ b=0.\bullet$$

\begin{crl} For $N=N(Z(A)),$
$$Z^{left}(A)\subseteq Z^{r.ass}(N).$$
\end{crl}
{\bf Proof.} Evident.

\begin{crl} \label{askar1} 
For any $a_1,\ldots, a_{n-1}\in N(Z(A)),$ 
and $a_n\in A,$ $z\in Z(A),$ the following relation takes place
$$a_1\circ a_2\circ\cdots a_{n-1}\circ a_n\circ z=
(a_1\circ a_2\circ\cdots a_{n-1}\circ a_n)\circ z.$$
\end{crl}

{\bf Proof.} For $n=2,$ the statement follows from the lemma. 
Suppose that this is also true for $n-1.$ Then by our lemma
$$a_1\circ (a_2\circ \cdots (a_{n-1}\circ (a_n\circ z))\cdots )=$$
$$a_1\circ \{(a_2\circ \cdots (a_{n-1}\circ a_n)\cdots )\circ z\}=$$
$$\{a_1\circ (a_2\circ \cdots (a_{n-1}\circ a_n)\cdots )\}\circ z
.\bullet$$

\begin{prp} \label{as} Let $U$ be a right antisymmetric $A-$module and 
$$A\cup U\rightarrow A$$
be a pairing of $A-$modules:
\begin{equation}
a\circ(b\cup u)=(a\circ b)\cup u,\label{cup1}
\end{equation}
$$(a\cup u)\circ b= (a\circ b)\cup u+a\circ(u\circ b),$$
for any $a,b\in A, u\in U$ (about cup-products see~\cite{DzhumaKos}).
Suppose that any element of $A$ can be presented 
by a cup-product as a $z\cup u,$ for some $u\in U$ and $z\in Z(A).$ Then 
for any $a_1,\ldots, a_{n-1}\in N(Z(A)),$ 
and $a_n\in A,$ 
$$a_1\circ a_2\circ\cdots a_{n-1}\circ a_n\circ a=
(a_1\circ a_2\circ\cdots a_{n-1}\circ a_n)\circ a.$$
\end{prp}

{\bf Proof.} Let $a=z\cup u.$ Then by~(\ref{cup1}), and corollary~\ref{askar1},
$$a_1\circ (a_2\circ \cdots (a_{n-1}\circ (a_n\circ a))\cdots )=$$
$$a_1\circ \{(a_2\circ \cdots (a_{n-1}\circ ((a_n\circ z)\cup u))\}=$$
$$a_1\circ (a_2\circ \cdots (a_{n-1}\circ (a_n\circ z))\cdots )\cup u=$$
$$\{a_1\circ (a_2\circ \cdots (a_{n-1}\circ a_n)\cdots )\}\circ z)\cup u=$$
$$\{a_1\circ (a_2\circ \cdots (a_{n-1}\circ a_n)\cdots )\}\circ (z\cup u)=$$
$$\{a_1\circ (a_2\circ \cdots (a_{n-1}\circ a_n)\cdots )\}\circ a.\bullet$$

\begin{prp} $Z^{l.ass}(A)\subseteq N(Z(A)).$
\end{prp}
{\bf Proof.} Let $a\in Z^{l.ass}(A).$ Then for any $z\in Z(A),$
$$(a\circ z)\circ b=a\circ(z\circ b)-(a,b,z)=0,$$
for any $b\in A.\bullet$

{\bf Example 1.} Any associative algebra is right-symmetric.
As associative algebra the matrix algebras $gl_n$ gives us examples
of right-symmetric algebras. 

{\bf Example 2.}
Let us give associative commutative algebra $U$ with commuting derivations 
${\cal D}=\{\der_i, i=1,\ldots,n\}.$ Then an algebra of derivations $UD=
\{u\der_i: u\in U, \der_i\in {\cal D}\}$ with multiplication $u\der_i\circ
v\der_j=
v\der_j(u)\der_i,$ is right-symmetric.
Because of Lie algebras for $U{\cal D}$ are Witt algebras:
$$[u\der_i,v\der_j]=-u\der_i(v)\der_j+v\der_j(u)\der_i,$$
we call such algebras as right-symmetric Witt algebras.  

Let $\Gamma_n$ be a set of $n-$typles $\alpha=(\alpha_1,\ldots,\alpha_n),$ 
where $\alpha_i$ are integers. Let $\Gamma_n^+$ be its subset consisting of 
such $\alpha,$  that $\alpha_i\ge 0, i=1,\ldots n.$ In the case of 
$p=char {\cal K}>0,$ 
we consider a subset $\Gamma_n({\bf m})=\{\alpha: 0\le \alpha_i<p^{m_i}, 
i=1,\ldots,n\},$ where ${\bf m}=(m_1,\ldots, m_n), m_i>0, m_i\in {\bf Z}, 
i=1,\ldots,n.$

For $char\,{\cal K}=0$ suppose that
$$U={\cal K}[x^{\pm 1}_1,\ldots,x^{\pm 1}_n]=\{x^{\alpha}=\prod_{i=1}^kx_i^{\alpha_i}:
\alpha\in \Gamma_n\}$$
\noindent is an algebra of Laurent polynomials and
$$U^+={\cal K}[x_1,\ldots,x_n]=\{x^{\alpha}: \alpha\in \Gamma_n^+\},$$
its subalgebra of polynomials. 

Let
$$O_n({\bf m})=\{x^{(\alpha)}=\prod_{i}x_i^{(\alpha_i)}:
\alpha\in \Gamma_n({\bf m}), i=1,\ldots,n\}$$

\noindent be divided power algebra if $char\,{\cal K}=p>0.$ Recall that  
$O_n({\bf m})$ is $p^m$~-dimensional and the muliplication is given by
$$x^{(\alpha)}x^{(\beta)}={{\alpha+\beta}\choose \alpha}
x^{(\alpha+\beta)},$$

\noindent where $m=\sum_im_i,$ and
$${{\alpha+\beta}\choose \alpha}=\prod_{i}{{\alpha_i+\beta_i}\choose \alpha_i},
\; {n \choose l}= {n!\over {l!(n-l)!}}, \; n,l\in {\bf Z}_+.$$

Let $\eps_i=(0,\ldots,\mathop{1}\limits_{i},\ldots,0).$ Define
$\der_i$ as a derivation of $U,$
$$\der_i(x^\alpha)=\alpha_ix^{\alpha-\eps_i}, \; p=0,$$
$$\der_i(x^{(\alpha)})=x^{(\alpha-\eps_i)}, \; p>0.$$

Denote right-symmetric algebras $U{\cal D}, U^+{\cal D}$ for 
$U={\cal K}[x^{\pm 1},\ldots ,x_n^{\pm 1}]$ as $W_n^{rsym}$ and $W_n^{+rsym}.$  
Denote, similarly, right-symmetric algebra $O_n({\bf m}){\cal D}$ as 
$W_n({\bf m})^{rsym}.$ 
Notice that Lie algebra $W_n$ is isomorphic to a Lie algebra of formal 
vector fields on $n-$dimensional torus and $W_n^+$ is isomorphic to a Lie 
algebra of formal vector fields on ${\cal K}^n.$ As in the case of Lie algebras,
$A=W_n^{rsym}$ has a grading 
$$A=\oplus_{k}A_k,\;\; A_k\circ A_l\in A_{k+l}, \;\; k,l\in {\bf Z},$$
$$A_k=\{x^{\alpha}\der_i:
|x^{(\alpha})|=|\alpha|=\sum_{i=1}^n\alpha_i=k+1\}.$$
This grading induces gradings in $W_n^{+rsym}$ and $W_n^{rsym}({\bf m}).$ 

{\bf Example 3.} Let $A$ be an associative algebra, $C^*(A,A)=\oplus_k C^k(A,A),$ 
and $C^k(A,A)=\{\psi: A\times \cdots\times A\rightarrow A\}$ be a space of
polylinear maps with $k-$arguments, if $k>0,$ $C^0(A,A)=A,$ and $C^k(A,A)=0,$  
if $k<0.$ 
Endow $C^*(A,A)$ by a "shuffle-product" multiplication:
$$C^*(A,A)\times C^*(A,A)\rightarrow C^{*-1}(A,A),$$
$$\psi\circ \phi(a_1,\ldots,a_{k+l+1})=$$
$$\sum_{\begin{array}{c}\sigma\in Sym_{k+l+1},\\
\sigma(1)<\cdots <\sigma(k+1),\\
\sigma(k+2)<\cdots <\sigma(k+l+1)
\end{array}}
\psi(\phi(a_{\sigma(1)},\ldots,a_{\sigma(k+1)}),a_{\sigma(k+2)},\ldots,a_{\sigma(k+l+1)}),$$
where
$\psi\in C^{k+1}(A,A), \phi\in C^{l+1}(A,A), \psi\circ \phi\in C^{k+l+1}(A,A), k,l\ge 0.$

Let $\Delta_k\in C^{k+1}(A,A), k\ge 0,$ be a standard skew-symmetric polynomial:
$$\Delta_k(a_1,\ldots,a_{k+1})=\sum_{\sigma\in Sym_{k+1}}
sign\,\sigma\,a_{\sigma(1)}\circ\cdots \circ a_{\sigma(k+1)}.$$

Then~\cite{Dzhuma88}
$$\Delta_i\circ \Delta_{2k}=(i+1)\Delta_{2k+i},$$
$$\Delta_{2k+1}\circ\Delta_{2l+1}=0,$$
$$\Delta_{2k}\circ \Delta_{2l+1}=\Delta_{2k+2l+1},$$
for any $k,l,i\ge 0.$

Therefore, the algebra of standard polynomials under shuffle-product is isomorphic
to the right-symmetric algebra $A_0\oplus A_1,$ such that 
$$A_0=\{e_i: i\ge 0\}, A_0\circ A_0\subseteq A_0,$$
$$ e_i\circ e_j=(i+1/2)e_{i+j}, 0\le i,j,$$
$$A_1=\{x^{j+1}: j\ge 0\}, A_1\circ A_1=0, A_0\circ A_1\subseteq A_1,
A_1\circ A_0\subseteq A_1,$$
$$x^{i+1}\circ x^{j+1}=0, e_i\circ x^{j+1}=(1/2)x^{i+j+1},\;\; i,j\ge 0.$$ 
The isomorphism can be given by
$$e_i\mapsto \Delta_{2i}/2, \;\; x^{j+1}\mapsto \Delta_{2j+1},$$
where $i,j=0,1,2,\ldots .$

This algebra has also  multiplication $\cup: C^*(A,A)\otimes C^*(A,A)
\rightarrow C^*(A,A),$ called as  a cup-product:
$$\psi\cup\phi(a_1,\ldots,a_{k+l})=\sum_{\begin{array}{c}
\sigma\in Sym_{k+l}\\
\sigma(1)<\cdots \sigma(k),\\
\sigma(k+1)<\cdots \sigma(k+l)
\end{array}}
\psi(a_{\sigma(1)},\ldots,a_{\sigma(k)})
\phi(a_{\sigma(k+1)},\ldots,a_{\sigma(k+l)}).
$$
Then
$$\Delta_i\cup\Delta_j=\Delta_{i+j+1},$$
for any $i,j\ge 0.$ In particular, a subalgebra generated by 
$\Delta_i,$ is a commutative associative algebra under cup-product.
These multiplications satisfy the following conditions
$$(\Delta_i\cup \Delta_j)\circ \Delta_k= (-1)^{k(j-1)}
(\Delta_i\circ \Delta_k)\cup \Delta_j+\Delta_i\cup(\Delta_j\circ \Delta_k),$$
$$(a\cup b)\cup c=a\cup(b\cup c),$$ 
$$a\cup b=b\cup a,$$
where $b\in C^k(A,A), c\in C^j(A,A).$ So, an algebra of standard polynomials has a
structure of Poissson-Novikov algebras in the sense of \cite{DzhumaComptes}.

\section{Novikov algebras}
A right-symmetric algebra $A$ is called (right) Novikov, if it satisfies
the following identity
$$a\circ (b\circ c)=b\circ(a\circ c), \;\; \forall a,b,c\in A.$$
Notice that algebras  $W_1^{rsym}, W_1^{+rsym}, W_1^{rsym}({m})$ are 
Novikov. The natural question about generalisation of Novikov identities that 
includes the case $n>1$ is arises. We suggest two ways to solve  
this problem.

In the first way, we will consider two kinds of multiplications on Witt algebras
$$u\der_i\circ v\der_j=v\der_j(u)\der_i,$$
$$u\der_i\ast v\der_j=\der_i(u)v\der_j.$$
They satisfy the following identities
$$a\circ (b\circ c)-(a\circ b)\circ c-a\circ(c\circ b)+(a\circ c)\circ b=0,$$
$$a\ast (b\ast c)-b\ast(a\ast c)=0,$$
$$a\circ (b\ast c)-b\ast(a\circ c)=0,$$
$$(a\ast b-b\ast a- a\circ b+b\circ a)\ast c=0,$$
$$(a\circ b- b\circ a)\ast c+a\ast(c\circ b)-(a\ast c)\circ b
-b\ast(c\circ a)+(b\ast c)\circ a=0.$$
In the case of $n=1$ these multiplications coincide and all identities are 
reduced to two: right-symmetric identity and Novikov identity. 
So, an algebra $A$ with multiplications $\circ$ and $\ast$ can be considered
as a Novikov algebras in the general case.

The second way concerns identities of right-symmetric Witt algebras. 

\section{Right-symmetric identities}
Let $R_n={\cal K}<t_1,\ldots,t_n>$ be a free algebra in the category of 
right-symmetric algebras generated on the set $\{t_1,\ldots,t_n\}.$ 
It can be defined in the following way. Consider a magma 
algebra $M<X_1,\ldots,X_n>$ generated on the set $\{X_1,\ldots,X_n\}$
~\cite{Bourbaki}.
Then 
$R_n$ is a factor algebra of $M<X_1,\ldots,X_n>$ over an ideal generated 
by associators $(X_i,X_j,X_s),$ $i,j,s=1,\ldots, n.$ Elements of $R_n$ 
corresponding to $X_1,\ldots,X_n$ are denoted by $t_1,\ldots, t_n.$ 

A polynomial $f(t_1,\ldots, t_n)\in {\cal K}<t_1,\ldots,t_n>$ is called 
{\it a right-symmetric identity} on $A,$ if $f(a_1,\ldots,a_n)=0,$ for any 
$a_1,\ldots,a_n\in A.$ If $f(a_1,\ldots,a_n)=0,$ is an identity on $A,$ then 
$f(a_{\sigma(1)},\ldots,a_{\sigma(n)})=0$ is also an identity on $A,$ for any
permutation $\sigma\in Sym_k.$ 

Further, denotions like $x_1\circ x_2\circ \cdots \circ x_{k-1}\circ x_k$ 
will mean an element $x_1\circ (x_2\circ (\cdots \circ(x_{k-1}\circ x_k)
\cdots)).$ For any finite sequence of integers ${\bf i}=(i_1,i_2,\ldots,i_k),$ 
with $i_l=1,\ldots,n,$ and $l=1,\ldots, k,$ where $k$ is any integer, 
set $t^{\bf i}=t_{i_1}\circ \ldots t_{i_k}.$ Elements of the form 
$\lambda_{\bf i}t^{\bf i},$ where $\lambda_{\bf i}\in {\cal K},$ 
are called as (left) monoms. If $\lambda_{\bf i}=0,$ then it is a trivial 
monom. If we would like to pay attention to monoms with $\lambda_{\bf i}\ne 0,$ 
then we call $\lambda_{\bf i}t^{\bf i}$ as a nontrivial monom. The sum of 
monoms is called as a (left) polynomial. A space of (left) polynomials is 
denoted by $R_n^{left}.$ 
In the universal enveloping algebra $U(R_n)$ they 
correspond to polynomials generated by $l_{t_i}, t_i\in R_n.$ {\it A degree} of 
the nontrivial monom $\lambda_{\bf i}t_{i_1}\circ \cdots\circ t_{i_k},$ by definition 
is $k:$
$$deg\,t^{\bf i}= |{\bf i}|=\sum_{l} i_l.$$
For a polynomial 
$$f=\sum_{\bf i} \lambda_{\bf i}t^{\bf i}\in R^{left}_n,$$ 
we will say that $f$ has a monom $\lambda_{\bf i}t^{\bf i},$ 
if $\lambda_{\bf i}\ne 0.$ The maximal degree of the nontrivial 
monoms of $f$ is called as a {\it degree} of the polynomial $f:$
$$def\,f= max\{|{\bf i}| : \lambda_{\bf i}\ne 0\}.$$

Suppose that  $A$ is a graded right-symmetric algebra:
$$A=\oplus_iA_i, \;\;\; \; A_i\circ A_j\subseteq A_{i+j}.$$
Elements of $A_i$ are called as homogeneous elements of the weight $i.$
Notation $|a|=i,$ will means that $a$ is a homogeneous element of $A$ and 
$a\in A_i.$ The element obtained from $t^{\bf i}$ by substituting 
$t_i:=a_i,$ we denote as $a^{\bf i}.$ 
For graded right symmetric algebra $A$ and 
for any monom $t^{\bf i},$ where ${\bf i}=(i_1,\ldots,i_k),$  it is evident that 
\begin{equation}
a^{\bf i}\in A_{\sum |a_{i_1}|+\cdots +|a_{i_k}|} \label{grading},
\end{equation}

In our paper we deal only in left polynomial identities.

The following left polynomial of $R^{left}_{k+1}$ is called 
a {\it right-symmetric standard polynomial} of degree $k+1:$
$$s_k^{rsym}(t_1,\ldots,t_{k+1})=\sum_{\sigma\in Sym_{k}}
sign\,\sigma t_{\sigma(1)}\circ\ldots \circ t_{\sigma(k)}\circ t_{k+1}.$$
Notice that it can be considered as an associative standard polynomial
of degree $k$ for left multiplications in the universal enveloping
algebra $U(R_{k+1})$~\cite{DzhumaKos}. If
$$s_k^{ass}(l_{t_1},\ldots,l_{t_k})=\sum_{\sigma\in Sym_k}sign\,\sigma\,
l_{t_{\sigma(1)}}\ldots l_{t_{\sigma(k)}},$$
then
$$s_k^{rsym}(t_1,\ldots,t_k,t_{k+1})=(t_{k+1})s_k^{ass}(l_{t_1},\ldots,l_{t_k}).$$

If $f(a_1,\ldots,a_k)=0$ is a right-symmetric polynomial identity of $A,$ then
the polynomial $g$ obtained from $f$ by permutation of arguments
$g(t_1,\ldots,t_k)=f(t_{\sigma(1)},\ldots,t_{\sigma(k)}), \;\; \sigma\in Sym_k,$
give also a right-symmetric polynomial identity. Denote $g$ as $\sigma\,f.$
Let $\tau_l\in Sym_{2n+1},\; l=0,1,\ldots, 2n,$ be the cyclic permutation
of elements ${0,1,\ldots,2n}:$ $\tau_l=(l,0,1,\ldots,l-1).$

The following is our main result.

\begin{th}\label{Askar} Let $A$ be one of the following right-symmetric algebras
$W_n^{rsym} (p=0), W^{+rsym}_n (p=0), W^{rsym}_n({\bf m}) (p>0).$

i) $A$ satisfies the right-symmetric standard identity of degree $2n+1:$
$$\sum_{\sigma\in Sym{2n}}sign\,{\sigma}\,a_{\sigma(1)}\circ
a_{\sigma(2)}\circ\cdots a_{\sigma(2n)}\circ a_{2n+1}=0,$$
for any $a_1,\ldots,a_{2n+1}\in A.$
In particular, the polynoms $\tau_ls_{2n}^{rsym}, l=0,1,\ldots,2n,$   
are also right-symmetric identities of $A.$

ii) If $char\,{\cal K}=0,$ then the algebra $A=W_n^{rsym},$ or $W_n^{+rsym},$
has no left polynomial identity of degree less that $2n+1.$ If $p>0,$ then
the algebra $W_n({\bf m})$ has no multilinear polynomial identity of degree
less than $2n+1.$

iii) Let $f$ be any left polynomial identity of $A$ of degree $2n+1.$
In the case of $W_n({\bf m}), p>0,$ we suppose that $f$ is also multilinear.
Then $f$ is a linear combination of $\tau_ls_{2n}^{rsym},\; l=0,1,\ldots,2n.$   
\end{th}

\section{Identities of right-symmetric  Witt algebras.}
In this section $A=W_n^{rsym}, W_n^{+rsym},$ or $W_n^{rsym}({\bf m}), p>0.$ 
Define  $s_{k,r}\in
C^k(A,A), r=1,\ldots,n,$ by
$$s_{k,r}(u_1\der_{i_1},\ldots,u_k\der_{i_k}) =
\sum_{\sigma\in Sym_k}sign(\sigma)\der_r(u_{\sigma(1)})
\der_{i_{\sigma(1)}}(u_{\sigma(2)})\ldots 
\der_{i_{\sigma(k-1)}}(u_{\sigma(k)})\der_{i_{\sigma(k)}}.$$
Let 
$$A^+_{-1}=\{\der_i: i=1,\ldots,n\},$$
$$A^+_0=\{x_i\der_j: i,j,=1,\ldots,n\}.$$

\begin{lm}\label{Witt}\cite{DzhumaKos}
$$Z(A)=A^+_{-1},$$
$$N(Z(A))=A^+_{-1}\oplus A^+_{0}.$$
In particular, $N(Z(A))$ has a subalgebra $A^+_0$ isomorphic to  $Mat_n.$
It is an associative subalgebra of $A.$
\end{lm}

The proof of theorem~\ref{Askar} is based on the following observation.
\begin{th} \label{basic} 
For any $r=1,\ldots,n,$ and $a_1,\ldots,a_{2n}\in A,$ 
$$s_{2n,r}(a_1,\ldots,a_{2n})=0.\bullet$$
\end{th}

Let $U={\cal K}[x^{\pm 1}_1,\ldots,x_n^{\pm 1}], 
{\cal K}[x_1,\ldots,x_n],$ or $O_n({\bf m}),$ 
if $A=W_n, W_n^+,$ or $W_n({\bf m}).$
It is enough to prove the theorem for a basic elements $a_i.$ 
Recall that they have the form $u\der_i,$ where $u=x^{\alpha}, \alpha\in \Gamma_n.$ 
For the proof of our theorem we need some lemmas.

\begin{lm} \label{one} 
If one of the basic elements $a_1,\ldots,a_k$ belongs to $A^+_{-1},$ then 
$$s_{k,r}(a_1,\ldots,a_{k})=0.$$
\end{lm}
{\bf Proof.} Evident.

\begin{lm} \label{two} For any $u,v\in U,$  and $1\le l \le k,$
$$s_{k,r}(a_1,\ldots, a_{l-1},
uv\der_{i_l},a_{l+1},\ldots,a_k)=$$
$$vs_{k,r}(a_1,\ldots, a_{l-1},
u\der_{i_l},a_{l+1},\ldots,a_k)+
us_{k,r}(a_1,\ldots, a_{l-1}, v\der_{i_l},a_{l+1},\ldots,a_k).$$
\end{lm}
{\bf Proof.} Evident.

\begin{lm} \label{three} If 
$s_{k,r}(a_1,\ldots,a_k)=0,$ for all $a_1,\ldots,
a_k\in W^{+rsym}_n,$ then 
$s_{k,r}(a_1,\ldots,a_k)=0,$ 
for all $a_1,\ldots,
a_k\in W^{rsym}_n.$ 
\end{lm}

{\bf Proof.} Let $l$ be a number of elements of $W_n$ in the set 
$\{a_1,\ldots,a_k\}.$ We use induction reasoning on $l=n,\ldots, 0.$ 

By the assumption of our lemma the base of induction is true. 
Suppose that the statement is true for $l$ and 
$l-1$ elements of the arguments set ${a_1,\ldots,a_k}$ belongs to 
$W_n^{+rsym}.$ Since the polynom $s_{k,r}$ is skew-symmetric, 
we can assume that $a_1,\ldots,a_{l-1}\in W^{+rsym}_n, $ and $a_{l}$ 
has a form $uv\der_i,$ for some $u=x^{-\alpha}, v=x^{\beta}, 
\alpha, \beta\in \Gamma_n^+.$ 

Present $1$ as $uu^{-1}.$ 
According lemma~\ref{one} and lemma~\ref{two}
$$0=s_{k,r}(a_1,\ldots,a_{l-1},\der_{i_l},
\ldots,a_{k})=$$
$$us_{k,r}(a_1,\ldots,a_{l-1},u^{-1}\der_{i_l},
\ldots,a_{k})+
u^{-1}s_{k,r}(a_1,\ldots,a_{l-1},u\der_{i_l}, \ldots,a_{k}).$$
Since $u^{-1}=x^{\alpha}\in U^+,$ according the inductive hypothesis, \\
$s_{k,r}(a_1,\ldots,a_{l-1},u^{-1}\der_{i_l}, \ldots,a_{k})=0.$ 
Therefore,
$$u^{-1}s_{k,r}(a_1,\ldots,a_{l-1},u\der_{i_l}, \ldots,a_{k})=0,$$ 
and
$$s_{k,r}(a_1,\ldots,a_{l-1},u\der_{i_l}, \ldots,a_{k})=
uu^{-1}s_{k,r}(a_1,\ldots,a_{l-1},u\der_{i_l}, \ldots,a_{k})=0.$$

Using once again lemma~\ref{two} and the inductive hypothesis
we have
$$s_{k,r}(a_1,\ldots,a_{l-1},uv\der_{i_l},
\ldots,a_{k})=$$
$$us_{k,r}(a_1,\ldots,a_{l-1},v\der_{i_l}, \ldots,a_{k})+
vs_{k,r}(a_1,\ldots,a_{l-1},u\der_{i_l}, \ldots,a_{k})=0.$$

So, the induction transfer is possible and the statement is proved.

\begin{lm}\label{four} If $f\in R^{left}_k,$ satisfy the condition
$$f(a_1,\ldots,a_k)=0,$$
for any $a_1,\ldots,a_k\in A^+_{-1}\oplus A^+_{0},$ then
$$f(a_1,\ldots,a_k)=0,$$
for any $a_1,\ldots,a_k\in A.$
\end{lm}

{\bf Proof.}
By lemma~\ref{three} it is enough to consider the cases of 
$A=W_n^+,$ and $A=W_n({\bf m}).$ Our arguments in both of these 
cases are similar and for definitely we assume that $A=W_n^+.$ 

If $|a_l|< 0,$ for some $1\le l \le 2n,$  by the condition 
of our lemma the statement is true. Suppose that $|a_l|\ge 0,$ for all $l=1,\ldots,n.$ 
We will argue by induction on $q=|a_1|+\cdots+|a_{2n}|.$ 

For $q=0,$ we have $|a_1|=\cdots|a_{2n}|=0.$ According to the condition of 
our lemma the statement is true.   

Suppose that for $q-1\ge 0, $ the statement is true and $a_1,\ldots,a_{2n}\in W_n^+,$ 
such that 
$|a_1|+\cdots+|a_{2n}|=q, \;\; |a_1|\ge 0,\ldots, |a_{2n}|\ge 0.$ 
We have
$$\der_t(s_{2n,r}(a_1,\ldots,a_{2n}))=\sum_{l=1}^{2n} 
s_{k,r}(a_1,\ldots,a_{l-1}, \der_t(u_{l})\der_{i_l}, \ldots,a_{2n}),$$
and
$$s_{k,r}(a_1,\ldots,a_{l-1}, \der_t(u_{l})\der_{i_l}, \ldots,a_{2n})=0,$$
since
$$|a_1|+\cdots+|a_{l-1}+|\der_t(u_l)\der_{i_l}|+\cdots+|a_{2n}|=q-1,
\;\;\;
 |\der_t(u_l)\der_{i_l}|\ge 0,$$

\noindent or
$$|\der_t(u_l)\der_{i_l}|=-1.$$

\noindent So, for any $1\le t\le n,$
\begin{equation}
\der_t(s_{2n,r}(a_1,\ldots,a_{l-1}, \der_t(u_{l})\der_{i_l}, \ldots,a_{2n}))=0.
\label{equ}
\end{equation}

\noindent Because of
$$\der_t(w)=0, w\in U, i=1,\ldots,n\Rightarrow 
w= \lambda 1, \mbox{\; for some\;\;} \lambda\in {\cal K},$$
\noindent and
$$s_{2n,r}(a_1,\ldots,a_{2n})\in A_q, \; q>0,$$

\noindent from ~(\ref{equ}) we have,
$$s_{2n,r}(a_1,\ldots,a_{2n})=0.$$
Therefore, the induction transfer is possible and the theorem 
is proved completely.

{\bf Proof of theorem~\ref{basic}.} 
If one element of the set $\{a_1,\ldots,a_{2n}\}$ belongs to $A_{-1}^+,$
then by lemma~\ref{one}, 
$$s_{2n,r}(a_1,\ldots,a_{2n})=0.$$
Suppose now that $a_1,\ldots,a_{2n}\in A_0^+=\{x_i\der_j: i,j=1,\ldots,n\}
\cong Mat_n.$ Let $a_l=x_{j_l}\der_{i_l}, \; l=1,\ldots,2n.$
According to the Amitsur-Levitzki theorem for any $r=1,\ldots,n,$
$$s_{2n,r}(a_1,\ldots,a_{2n})=$$  
$$\der_r\{\sum_{\sigma\in Sym{2n}}
sign\,\sigma\, x_{j_{\sigma{(1)}}}\der_{i_{\sigma(1)}}(x_{j_{\sigma(2)}})\cdots
\der_{i_{\sigma{(2n-1)}}}(x_{j_{\sigma(2n)}})\der_{i_{\sigma{(2n)}}}\}=0.$$
By lemma~\ref{four}, theorem~\ref{basic} is established.

{\bf Proof of theorem~\ref{Askar}.} 

i) Let $a_l=u_l\der_{i_l}, \; l=1,2,\ldots,2n, \;
a_{0}=u\der_r.$ Then by theorem~\ref{basic}
$$\sum_{\sigma\in Sym_{2n}}sign\,{\sigma}\,a_{\sigma(2n)}\circ
a_{\sigma(2)}\circ
\cdots a_{\sigma(1)}\circ a_{0}=$$
$$\sum_{\sigma\in Sym_{2n}}u\der_r(u_1)\der_{i_1}(u_2)\cdots 
\der_{2n-1}(u_{2n})\der_{2n}=$$
$$us_{2n,r}(a_1,\ldots,a_{2n})=0.$$

ii) If $p=0,$ and $A$ has nontrivial left polynomial identity of degree $d,$ then it has
multilinear nontrivial left polynomial identity of degree $\le d.$ 
The proof of this statement is based on the linearisation method and it
does not depend on the condition of associativity of $A.$ 
Let $f(x_1,\ldots,x_d)$ be a multilinear  identity of $A$ of degree
$\le 2n.$ Then $f(x_1,\ldots,x_n)$ has a monom of the form $\lambda
x_dx_{d-1}\ldots x_1, \; \lambda\ne 0.$ As in the case of matrices, if
$d<2n+1,$ we can take $a_1=\der_1,
a_2=x_1\der_1,a_3=x_1\der_2, a_4=x_2\der_2, \ldots
,a_{d}=x_m\der_m,$ if $d=2m,$ and $a_{d}=x_{m-1}\der_{m},$ if $d=2m-1.$ 
Then 
$$a_d\circ a_{d-1}\circ \cdots \circ a_1=\der_m,$$
and for any permutation $\sigma\in Sym_d,$ $\sigma\ne id,$ 
$$a_{\sigma(d)}\circ \ldots \circ a_{\sigma(1)}=0.$$
Therefore,
$$f(a_1,\ldots,a_d)=\lambda \der_m\ne 0.$$
Contradiction. 

iii) Let $f$ be a left polynomial of degree $2n+1$ and $f$ depends on 
valuables $t_1,\ldots,t_s,$ i.e. $f\in R^{left}_{s}.$ 
Suppose that $f=0$ is an identity on $A.$ By lemma~\ref{four} it is 
equivalent to say that $f=0$ is an identity on $N(Z(A))= 
A^+_{-1}\oplus A^+_{0}.$

Let 
\begin{equation}\label{fl}
f(t_1,\ldots,t_{s})=\sum_k\sum_{i_1,\ldots,i_{k}=1,\ldots,s}
\lambda_{i_1,\ldots,i_{k-1},i_k}
t_{i_1}\circ \cdots t_{i_{k-1}}\circ t_{i_k}.
\end{equation}
We would like to prove that $s=2n+1.$ If $f$ is multilinear, it is evident.
Since in the case of $A=W_n({\bf m})^{rsym},$ we suppose that $f$ is multilinear,
below we take $A=W_n^{rsym}, W_n^{+rsym}, p=0.$

Let us prove that $s\ne 1.$ Suppose that it is false. Let 
$$\mu_k=\lambda_{\mathop{\underbrace{1,\ldots,1}}
\limits_{k \mbox{\; times}}},
\mbox{\;and\;\;} t_1^k=\mathop{\underbrace{t_1\circ\cdots\circ t_1}
\limits_{k \mbox{\; times}}}.$$
Then
$$f(t_1)=\sum_k\mu_kt_1^k.$$
Substitute $a_1=x_1^{2}\der_1$ instead of $t_1.$  Since
$$a_1^{\circ k}=2^{k-1}x^{k+1}\der_1,$$
we have
$$\sum_k\mu_k2^{k-1}x_1^{k+1}\der_1=0,$$
or
$$\mu_k2^{k-1}=0,$$
for any $k.$ Therefore, $\mu_k=0,$ for any $k,$ and $f$ is a trivial polynomial,
i.e., $f=0.$ So, the case $s=1$ is not possible.

Suppose that all left polynomials of degree $2n+1$ with $s-1$ valuables are 
trivial.  Define 
$f_s\in R^{left}_{s}\subseteq M<t_1,\ldots,t_{s}>$ by 
$$f_s(t_1,\ldots,t_{s})=
\sum_k\sum_{i_1,\ldots,i_{k-1}=1,\ldots,s-1}
\lambda_{i_1,\ldots,i_{k-1},s}
t_{i_1}\circ \cdots t_{i_{k-1}}\circ t_s.$$
Denote by $\bar{f}_s$ the sum of monoms $\lambda_{\bf i}t^{\bf i},$ 
such that all coordinates of ${\bf i}=(i_1,\ldots,i_k)$ differ from $s.$
So, $f_s$ has monoms $\lambda_{\bf i}t^{\bf i},$ such that the entrance of 
$t_s$ is exactly 1 and $t_s$ appears only at the end of monoms. As far 
as $\bar{f}_l\in R^{left}_{s},$ 
its monoms do not depend from $t_s.$ So $\bar{f}_l\in R^{left}_{s-1}.$ 

It is evident that $a_{i_1}\circ \cdots \circ a_{i_k}=0,$ if 
$a_{i_l}\in Z(A),$ for some $l<k.$ Thus for any $z\in Z(A),$ 
we have
\begin{equation}
0=f(a_1,\ldots,a_{s-1},z)=
f_s(a_1,\ldots,a_{s-1},z)+
\bar{f}_s(a_1,\ldots, a_{s-1}).
\label{one}
\end{equation}
By corollary~\ref{askar1}, for any $a_1,\ldots,a_{s-1}\in N(Z(A)),$
$$f_s(a_1,\ldots,a_{s-1},z)=
F_s(a_1,\ldots,a_{s-1})\circ z,
$$
for some $F_s\in R^{left}_{s-1}.$ Namely,
$$F_s(x_1,\ldots, x_{s-1})=
\sum_k\sum_{i_1,\ldots,i_{k-1}=1,\ldots,{s-1}}
\lambda_{i_1,\ldots,i_{k-1},s}
t_{i_1}\circ \cdots t_{i_{k-1}},$$
if $f_s$ has a form~(\ref{fl}). 

Thus, the condition (\ref{one}) can be written in the following way
\begin{equation}
F_s(a_1,\ldots,a_{s-1})\circ z+
\bar{f}_l(a_1,\ldots, a_{s-1})=0.
\label{two}
\end{equation}
The first summand of the left hand of (\ref{two}) depends on $z$ linearly 
and the second summand does not depend on z. Therefore, 
\begin{equation}
F_s(a_1,\ldots,a_{s-1})\circ z=0,
\label{three}
\end{equation}
\begin{equation}
\bar{f}_s(a_1,\ldots, a_{s-1})=0,
\label{four}
\end{equation}
for any $z\in Z(A)$ and $a_1,\ldots,a_{s-1}\in N(Z(A)).$ 

By  (\ref{four}) and lemma~\ref{four} the polynomial $\bar{f}_s$ 
gives us an identity of degree no more that $2n+1.$ 
If this degree is less than $2n+1,$ according (ii), $\bar{f}_s$ should be 
a trivial polynomial. If $deg\,\bar{f}_s=2n+1,$ then by the inductive 
suggestion $\bar{f}_s,$ as a polynomial 
with $s-1$ valuables, also should be trivial. 

Notice that $deg\,f_s\le 2n+1,$ and $def\,F_s\le 2n.$
Since the centraliser of $Z(A)$ coincides with $Z(A)$ (see lemma~\ref{Witt}) 
and $Z(A)\subset A_{-1},$ by (\ref{grading}) from (\ref{three}) we  obtain  
that
$$F_s(a_1,\ldots,a_{s-1})=0$$
is an identity for $N(Z(A)).$ In particular, it is an identity on $A^+_0\subset 
N(Z(A))$ (recall that $N(Z(A))=A^+_{-1}\oplus A^+_0$). We have mentioned 
that $A^+_0=\{x_i\der_j: i,j=1,\ldots,n\}\cong Mat_n.$ So, $F_s=0$ is an 
identity on $Mat_{n}$ of a degree no more than $2n.$

By the Amitsur-Levitzki theorem, the minimal degree of nontrivial polynomial 
identity on $Mat_n$ is $2n$ and any such polynomial should be a standard 
polynomial of degree $2n$ up to scalar. Therefore, $s-1=2n.$ Moreover, 
$$F_{2n+1}(a_1,\ldots,a_{2n})=
\gamma_{2n+1}s_{2n}^{ass}(a_1,\ldots,a_{2n}),$$
where
$$s_{2n}^{ass}(t_1,\ldots,t_{2n})=
\sum_{\sigma\in Sym_{2n}}
sign\,\sigma\,t_{\sigma(1)}\circ \cdots\circ t_{\sigma(2n)}),$$
and $\gamma_{2n+1}\in {\cal K}.$ 

So, for 
$$g_{2n+1}(a_1,\ldots,a_{2n},a_{2n+1})=$$
$$f_{2n+1}(a_1,\ldots, a_{2n},a_{2n+1})-\gamma_{2n+1}\,s_{2n}^{rsym}
(a_1,\ldots,a_{2n},a_{2n+1})$$
and for any $a_{2n+1}=u\der_r\in A,$ $a_1,\ldots,a_{2n}\in A^+_{-1}\oplus A^+_0,$  
by proposition~\ref{as} we have
$$g_{2n+1}(a_1,\ldots,a_{2n},a_{2n+1})=$$
$$f_{2n+1}(a_1,\ldots, a_{2n},a_{2n+1})-\mu\,s_{2n}^{rsym}
(a_1,\ldots,a_{2n},a_{2n+1})=$$
$$u\der_r\{F_{2n+1}(a_1,\ldots, a_{2n})\}-
\gamma_{2n+1}\,u\der_r\{s_{2n}^{ass}(a_1,\ldots,a_{2n})\}=$$
$$u\der_r\{F_{2n+1}(a_1,\ldots, a_{2n})-s_{2n}^{ass}(a_1,\ldots,a_{2n})\}=$$
$$0.$$

By lemma~\ref{four}, 
$$g_{2n+1}(a_1,\ldots,a_{2n+1})=0,$$
for any $a_1,\ldots,a_{2n+1}\in A.$ Thus,
$$f_s=\gamma_s s_{2n}^{rsym},$$
for $s=2n+1.$

For any $l=1,\ldots,2n,$ let $f_l$ be a sum of monoms $\lambda_{\bf i}t^{\bf i}$ of $f,$ 
such that ${\bf i}=(i_1,\ldots,i_k)$ ends by $l$ and $i_1,\ldots,l_{k-1}\ne l:$
$$f_l(t_1,\ldots,t_{2n+1})=
\sum_k\sum_{i_1,\ldots,i_{k-1}=1,\ldots\hat{l},\ldots,2n+1}
\lambda_{i_1,\ldots,i_{k-1},l}
t_{i_1}\circ \cdots t_{i_{k-1}}\circ t_l.$$
Recall that the denotion $\hat{x}$ means that $x$ is omitted. 
Let $\bar{f}_l$ be a sum of monoms $\lambda_{\bf i}t^{\bf i},$ 
such that all coordinates of ${\bf i}=(i_1,\ldots,i_k)$ differ from $l.$
Reorder valuables $t_1,\ldots,t_l,\ldots,t_{2k+1},$ 
in a such way, that the last valuable will be $t_l.$ 
Repeating the arguments given above we can see that polynomial $\bar{f}_l$ is trivial  
and $f_l$ is equal to a polynomial obtained from $\gamma_ls_{2n}^{rsym}$ 
by permutation of arguments, where $\gamma_l\in {\cal K}.$ 

Therefore, $f_{2n+1}$ is a linear combination of standard polynomials
$\tau_ls_{2n}^{rsym}, l=0,1,\ldots,2n.$ $\bullet$

\end{document}